\documentclass[12pt,reqno]{amsart}
\usepackage{mathrsfs}
\usepackage{amsgen}
\usepackage{amscd}

\usepackage{amsmath}
\usepackage{latexsym}
\usepackage{amsfonts}
\usepackage{amssymb}
\usepackage{amsthm}
\usepackage{graphicx,epsfig}
\usepackage{color}
\usepackage{amsthm,amssymb}
\usepackage{graphicx}
\usepackage{enumerate}
\usepackage{amssymb,amsmath,graphicx,amsfonts,euscript}
\usepackage{color}
\usepackage{mathrsfs}
\usepackage{ulem}
\usepackage[margin=3cm, a4paper]{geometry}
\DeclareGraphicsExtensions{.eps}
\usepackage{amssymb,amsmath,graphicx,amsfonts,euscript,mathrsfs}
\usepackage{hyperref}

\usepackage{epstopdf}
\newtheorem{theorem}{Theorem}[section]
\newtheorem{lemma}[theorem]{Lemma}

\theoremstyle{definition}

\theoremstyle{remark}

\def\R{\mathbb{R}}
\def\Z{\mathbb{Z}}
\def\T{\mathbb{T}}
\def\C{\mathbb{C}}

\newcommand{\fe}{\mathrm{e}}

\newcommand{\bR}{{\mathbb R}}

\newcommand{\bT}{{\mathbb T}}

\newcommand{\bxi}{\boldsymbol\xi}

\newcommand{\minus}{\scalebox{0.75}[1.0]{$-$}}

\numberwithin{equation}{section}

\begin{document}

\title[Low-regularity integrator for NLS equation]{A second order low-regularity integrator for the nonlinear Schr\"odinger equation}

\author[A.~Ostermann]{Alexander Ostermann}
\address{\hspace*{-12pt}A.~Ostermann: Department of Mathematics, University of Innsbruck, Technikerstra\ss e 13, 6020 Innsbruck, Austria}
\email{alexander.ostermann@uibk.ac.at}

\author[Y.~Wu]{Yifei Wu}
\address{\hspace*{-12pt}Y.~Wu: Center for Applied Mathematics, Tianjin University, 300072, Tianjin, China}
\email{yerfmath@gmail.com}

\author[F. Yao]{Fangyan Yao}
\address{\hspace*{-12pt}F.~Yao: School of Mathematical Sciences,
South China University of Technology, Guangzhou, Guangdong 510640, P. R. China}
\email{yfy1357@126.com}


\date{}

\dedicatory{}

\begin{abstract}
In this paper, we analyse a new exponential-type integrator for the nonlinear cubic Schr\"odinger equation on the $d$ dimensional torus $\T^d$. The scheme has recently also been derived in a wider context of decorated trees in \cite{bs}. It is explicit and efficient to implement. Here, we present an alternative derivation, and we give a rigorous error analysis. In particular, we prove second-order convergence in $H^\gamma(\T^d)$ for initial data in $H^{\gamma+2}(\T^d)$ for any $\gamma > d/2$. This improves the previous work in~\cite{lownls2}.

The design of the scheme is based on a new method to approximate the nonlinear frequency interaction. This allows us to deal with the complex resonance structure in arbitrary dimensions. Numerical experiments that are in line with the theoretical result complement this work.
\end{abstract}

\maketitle

\section{Introduction}\label{sec:introduction}
The nonlinear Schr\"{o}dinger equation (NLS) arises as a model equation in several areas of physics, see, e.g., Sulem and Sulem \cite{sulem}. In this paper, we are concerned with the numerical integration of the NLS equation on a $d$ dimensional torus:
\begin{equation}\label{model}
 \left\{\begin{aligned}
& i\partial_tu(t,x)+\Delta u(t,x)
 +\lambda|u(t,x)|^2u(t,x)=0,
 \quad t>0,\ x\in\bT^d,\\
 &u(0,x)=u_0(x),\quad x\in\bT^d,
 \end{aligned}\right.
\end{equation}
where $\bT=(0,2\pi)$, $\lambda=\pm1$, $u=u(t,x):\bR^{+}\times\bT^d\to\C$ is the sought-after solution, and $u_0\in H^\gamma(\bT^d)$ for some $\gamma\ge 0$ is the given initial data. Here we only consider the case $\lambda=1$; the case $\lambda=-1$ can be treated in exactly the same way. Note that the well-posedness of the nonlinear Schr\"{o}dinger equation in $H^\gamma(\T^d)$ has been established for $\gamma>\frac d2-1$. For details, we refer to \cite{Bo}.

Numerical aspects of the NLS equation have been studied by many authors. A considerable amount of literature has been published on splitting methods and exponential integration methods. For a general introduction to these methods, we refer to \cite{faou,hairer,Hochbruck,holden,Splitting}. It is well known that schemes of arbitrarily high order can be constructed by assuming that the solution of \eqref{model} is smooth enough. For instance, second order convergence in $H^\gamma$ was obtained by requiring four additional derivatives of the solution for the Strang splitting scheme in \cite{Lubich}. Further convergence results for semilinear Schr\"{o}dinger equations can be found, e.g., in \cite{besse,cano,celledoni,cohen,djuardin,ignat,jahnke,thalhammer}.

For classical methods and their analysis, strong regularity assumptions are unavoidable. Recently, however, so-called low-regularity integrators have emerged as a powerful tool for reducing the regularity requirements. The first breakthrough was made in \cite{lownls}, where the authors introduced a new exponential-type numerical scheme and achieved first-order convergence in $H^\gamma(\T^d)$ for $H^{\gamma+1}(\T^d)$ initial data. Later, a first-order integrator was proposed in \cite{wu20f}. It converges in $H^\gamma(\T)$ without any loss of regularity and conserves mass up to order five. A second-order Fourier-type integrator was given by Kn\"{o}ller, Ostermann and Schratz \cite{lownls2}. The integrator is based on the variation-of-constants formula and makes use of certain resonance based approximations in Fourier space. For second-order convergence, the scheme requires two additional derivatives of the solution in one space dimension and three derivatives in higher space dimensions. In this paper, we present and analyse an improved integrator which enables us to get the desired second-order accuracy with only two additional bounded spatial derivatives in dimensions $d\ge1$.

There are two main difficulties in designing low-regularity integrators. The first one is to control the spatial derivatives in the approximation while keeping the nonlinearity point-wise defined in physical space rather than in Fourier space. The second one is to overcome the difficulties caused by the complicated structure of resonances in higher dimensions. To explain this, let
$$
\bxi=(\xi^1, \cdots, \xi^d)\in \mathbb{Z}^d,\quad \bxi\cdot\boldsymbol\eta=\xi^1\eta^1+\cdots+ \xi^d\eta^d,\quad |\bxi|^2=\bxi\cdot\bxi.
$$
and consider the phase function
$$
\phi_3=|\bxi|^2+|\bxi_1|^2-|\bxi_2|^2-|\bxi_3|^2.
$$
In \cite{lownls2}, letting
\begin{equation*}
\alpha=2|\bxi_1|^2,\quad \beta=2\bxi_1\cdot\bxi_2+2\bxi_1\cdot\bxi_3+2\bxi_2\cdot\bxi_3,
\end{equation*}
the authors approximated the phase function by
\begin{align}\label{pre-approx}
\fe^{is\phi_3} =\fe^{is\alpha+is\beta} = \fe^{is\alpha}+\fe^{is\beta}-1+\mathcal R_1(\alpha,\beta,s),
\end{align}
where $|\mathcal R_1(\alpha,\beta,s)|\lesssim s^2|\alpha||\beta|$. This choice requires three additional derivatives in higher space dimensions for second-order convergence.

Now we explain our present approach, for which we consider a slightly more general situation. Assume that $\alpha$ has a ``good'' structure which means $\int_0^\tau\fe^{is\alpha}\,ds$  is point-wise defined (as in the above example) while $\beta$ has a ``bad'' structure but still has a low upper-bound, e.g., consisting of mixed derivatives (as in the example above). Then we employ the following approximation
\begin{align}\label{app}
\fe^{is(\alpha+\beta)}\approx\fe^{is\alpha}+i\beta \fe^{is\beta} \mathcal{M}_\tau\big(\fe^{i\alpha\cdot}\big),
\end{align}
where the operator $\mathcal{M}_\tau$ is defined by
\begin{align}\label{def:M}
\mathcal{M}_\tau(g)=\frac1\tau\int_0^\tau \sigma g(\sigma)\,d\sigma.
\end{align}
The mean $\mathcal{M}_\tau\big(\fe^{i\alpha\cdot}\big)$ can be regarded as an approximation to $\tau \fe^{i\alpha \tau}$ for small $\tau$.

Using the approximation \eqref{app} has several advantages. First, its integral has a point-wise interpretation in physical space,
\begin{align}\label{pointwise}
\int_0^\tau\Big(\fe^{is\alpha} +i\beta \fe^{is\beta} \mathcal{M}_\tau\big(\fe^{i\alpha\cdot}\big)\Big)\,ds
=\tau \varphi(i\tau \alpha)-\tau\big(\fe^{i\tau\beta}-1\big)\psi(i\tau \alpha),
\end{align}
see Lemma \ref{lem:idea-1} below. Here, the functions $\varphi$ and $\psi$ are defined as
\begin{align}\label{11}
\varphi(z)=\left\{\begin{aligned}
&\frac{e^z-1}{z},&\quad z\ne 0,\\
& \ 1,&\quad z= 0,
\end{aligned} \right. \qquad \quad
\psi(z)=\left\{\begin{aligned}
&\frac{e^z-1-ze^z}{z^2},&\quad z\ne 0,\\
& -\frac12,&\quad z= 0.
\end{aligned} \right.
\end{align}

Secondly, it is a high-order approximation and requires less regularity. Indeed, we get
\begin{align}\label{new-approx}
\int_0^\tau\fe^{is(\alpha+\beta)}\,ds=\tau \varphi(i\tau \alpha)-\tau\big(\fe^{i\tau\beta}-1\big)\psi(i\tau \alpha)+\mathcal R_2(\alpha,\beta,\tau),
\end{align}
where
$|\mathcal R_2(\alpha,\beta,\tau)|\lesssim \tau^3|\beta|^2$. This will be proved in Lemma \ref{lem:idea-1} below.
Relying on this structure, the scheme requires only two additional derivatives for $\tau^2$, which gives convergence in $H^\gamma(\T^d)$ for initial data in $H^{\gamma+2}(\T^d)$.

Finally, it does not require any specific structure of $\beta$. In particular, $\beta^{-1}$ is not contained in the expression \eqref{pointwise}. This is another advantage compared to \eqref{pre-approx}, for which the integration (or a further approximation) of $\int_0^\tau\fe^{is\beta}\,ds$ is needed.

Now we state the main result of this paper. We define the new low-regularity integrators with second-order accuracy as
\begin{equation}\label{LRI-discrete-1}
\begin{alignedat}{2}
u^0&=u_0,\\
u^{n+1}&=\fe^{i\tau\Delta}u^n&&+i\tau\fe^{i\tau\Delta}\Big\{\big[\varphi(-2i\tau \Delta)+\psi(-2i\tau \Delta)\big]\bar u^n \cdot(u^n)^2\Big\}\\
&&&-i\tau\big[\fe^{i\tau\Delta}\psi(-2i\tau \Delta)\bar u^n\big]\cdot(\fe^{i\tau\Delta}u^n)^2-\frac{\tau^2}{2}\fe^{i\tau\Delta}\Big[\big|u^n\big|^4 u^n\Big]
\end{alignedat}
\end{equation}
for $n\ge 0$. For this method, we have the following convergence result.
\begin{theorem}\label{thm:convergence}
Let $u^n$ be the numerical solution \eqref{LRI-discrete-1} of the Schr\"{o}dinger equation (\ref{model}) up to some fixed time $T>0$. Under the assumption that $u_0\in H^{\gamma+2}(\bT^d)$ for some $\gamma>\frac{d}{2}$, there exist constants $\tau_0,\,C>0$ such that for any $0<\tau\leq\tau_0$ it holds
\begin{equation}\label{main result}
\|u(t_n,\cdot)-u^n\|_{H^\gamma}\leq C\tau^2,\qquad 0\le n\tau \le T.
\end{equation}
The constants $\tau_0$ and $C$ only depend on $T$ and $\|u\|_{L^\infty((0,T);H^{\gamma+2}(\bT^d))}$.
\end{theorem}

Having finished the analysis of this paper, we became aware of the recent work \cite{bs} by Bruned and Schratz, in which low-regularity integrators for dispersive equations are discussed in a wider context. In particular, using the formalism of decorated trees, various numerical methods for the nonlinear Schr\"{o}dinger equation are proposed. The above method \eqref{LRI-discrete-1} is stated there in formula (5.17). Nevertheless, we give here alternative (and brief) derivation of the method because the employed approximations form the basis of our rigorous error analysis.

The paper is organized as follows.
In Section~\ref{sec:proof}, we introduce some notation and collect some useful lemmas. In Section~\ref{sec:construction}, we discuss the construction of the method and analyse accuracy and regularity requirements of each single approximation step. Collecting all these results, we proof our convergence result (Theorem~\ref{thm:convergence}) in Section~\ref{sec:ord-proof}. This theoretical result is illustrated with some numerical experiments in Section~\ref{sec:numerical}.

\section{Preliminaries}\label{sec:proof}

In this section, we introduce some notation, recall a result from harmonic analysis and give some elementary estimates. All this will be frequently used in the following sections.
\subsection{Some notation}\label{subsec1}
We start with notation, some of it borrowed from \cite{CKSTT-04-KDV}. We write $A\lesssim B$ or $B\gtrsim A$ to denote the statement that $A\leq CB$ for some constant $C>0$ which may vary from line to line but is independent of $\tau$ or $n$, and we write $A\sim B$ for $A\lesssim B\lesssim A$. We further denote
$$
\left<\bxi\right> = \sqrt{1+\bxi\cdot\bxi}, \qquad \bxi=(\xi^1, \cdots, \xi^d)\in \mathbb{Z}^d
$$
and define $(d\bxi)$ to be the normalized counting measure on $\Z^d$ such that
\begin{equation*}
\displaystyle\int
a(\bxi)\,(d\bxi)
=
\sum\limits_{\bxi\in
\Z^d} a(\bxi).
\end{equation*}
The Fourier transform of a function $f$ on $\T^d$ is defined by
$$
\hat{f}(\bxi)=\frac1{(2\pi)^d}\displaystyle\int_{\bT^d} \fe^{- i  \boldsymbol x\cdot\bxi}f(\boldsymbol x)\,d\boldsymbol x.
$$
Instead of $\hat f$, we
sometimes also write $\mathcal Ff$ or $\mathcal F(f)$.  The Fourier inversion formula takes the form
$$
f(\boldsymbol x)=\displaystyle\int \fe^{i  \boldsymbol x\cdot\bxi} \hat{f}(\bxi)\,(d\bxi).
$$
We recall the following  properties of the Fourier transform:
\begin{eqnarray*}
 &\|f\|_{L^2(\bT^d)}
 =(2\pi)^{\frac d2} \big\|\hat{f}\big\|_{L^2((d \bxi))} \quad \mbox{(Plancherel)}; \\
 &\displaystyle\langle f,g\rangle=\int_{\bT^d} f(\boldsymbol x)\overline{g(\boldsymbol x)}\,d\boldsymbol x
 =(2\pi)^d \displaystyle\int \hat{f}(\bxi)\overline{\hat{g}(\bxi)}\,(d\bxi)\quad \mbox{(Parseval)} ; \\
 & \widehat{(fg)}(\bxi)=\displaystyle\int
  \hat{f}(\bxi-\boldsymbol\eta)\hat{g}(\boldsymbol\eta) \,(d\boldsymbol\eta) \quad \mbox{(convolution)}.
\end{eqnarray*}
For the Sobolev space $H^s(\bT^d)$, $s\geq0$, we consider the equivalent norm
$$
\big\|f\big\|_{H^s(\bT^d)}=
\big\|J^sf\big\|_{L^2(\bT^d)}=(2\pi)^{\frac{d}{2} }\left\|(1+|\bxi|^2)^{\frac{s}{2}}\hat{f}(\bxi)\right\|_{L^2((d\bxi))},
$$
where $J^s=(1-\Delta)^\frac s2$.
\subsection{Some estimates}\label{subsec3}
First, we recall the following inequality, which was originally proved in \cite{Kato-Ponce}.
\begin{lemma}\label{lem:kato-Ponce}(Kato-Ponce inequality, \cite{Kato-Ponce}) The following inequalities hold:
\begin{itemize}
  \item[(i)]
  For any $ \gamma>\frac d2$ and $f,g\in H^{\gamma}$, we have
\begin{align*}
\|J^\gamma (fg)\|_{L^2}\lesssim \|f\|_{H^\gamma}\|g\|_{H^{\gamma}}.
\end{align*}
  \item[(ii)]
For any $\delta\ge 0, \gamma>\frac d2$ and $f\in H^{\delta+\gamma}$, $g\in H^{\delta}$, we have
\begin{align*}
\|J^\delta (fg)\|_{L^2}\lesssim \|f\|_{H^{\delta+\gamma}}\|g\|_{H^{\delta}}.
\end{align*}
\end{itemize}
\end{lemma}

The next lemma plays a crucial role in the analysis of this paper.
\begin{lemma}\label{lem:idea-1}
Let $\mathcal{M}_\tau$ be the operator defined in \eqref{def:M} and $\alpha,\beta\in \R$. Then, the following properties hold.
\begin{itemize}
\item[(i)] For $\varphi, \psi$ defined as in \eqref{11}, we have
\begin{align}\label{point-wise}
\int_0^\tau\Big(\fe^{is\alpha}
           +i\beta \fe^{is\beta} \mathcal{M}_\tau(\fe^{i\alpha\cdot})\Big)\,ds
=\tau \varphi(i\tau \alpha)-\tau\big(\fe^{i\tau\beta}-1\big)\psi(i\tau \alpha).
\end{align}
\item[(ii)] There exists a function $\mathcal R_2(\alpha,\beta,\tau)$ such that
\begin{align}\label{idea}
\int_0^\tau\fe^{is(\alpha+\beta)}\,ds=\int_0^\tau\Big(\fe^{is\alpha}
           +i\beta \fe^{is\beta} \mathcal{M}_\tau(\fe^{i\alpha\cdot})\Big)\,ds+\mathcal R_2(\alpha,\beta,\tau),
\end{align}
with
$|\mathcal R_2(\alpha,\beta,\tau)|\lesssim \tau^3|\beta|^2$.
\end{itemize}
\end{lemma}

\begin{proof}
(i)
The left-hand side of \eqref{point-wise} is equal to
\begin{align}
\int_0^\tau\fe^{is\alpha}\,ds
    +\int_0^\tau i\beta \fe^{is\beta}\mathcal{M}_\tau(\fe^{i\alpha\cdot})\,ds.\label{4.09-0424}
\end{align}
For the first term in \eqref{4.09-0424}, we have that
\begin{equation}\label{J1}
\int_0^\tau \fe^{is\alpha}\,ds =
\left.\begin{cases}
\dfrac{\fe^{i\tau\alpha}-1}{i\alpha},&\alpha\ne 0\\
\ \tau,&\alpha= 0\end{cases}
\right\}
= \tau \varphi(i\tau \alpha).
\end{equation}
For the second term in \eqref{4.09-0424}, using integration-by-parts, we find that
 \begin{equation}
\mathcal{M}_\tau(\fe^{i\alpha\cdot})= \frac{1}{\tau}\int_0^\tau \sigma\fe^{i\sigma\alpha}ds
 =\left\{\begin{aligned}
& \frac{\fe^{i\tau \alpha}}{i\alpha}
           +\frac{\fe^{i\tau \alpha}-1}{\tau \alpha^2},&\quad \alpha\ne 0,\\
& \ \tfrac12\tau,&\quad \alpha= 0.
\end{aligned} \right.
\end{equation}
Thus, from the definition of $\psi$, we infer that
\begin{align}\label{G-psi}
 \mathcal{M}_\tau(\fe^{i\alpha\cdot})=-\tau\psi(i\tau \alpha).
\end{align}
In addition,
\begin{align}\label{J2-2}
\int_0^\tau i\beta\fe^{is\beta}\,ds=\fe^{i\tau\beta}-1.
\end{align}
Therefore, combining  \eqref{J1} with \eqref{G-psi} and \eqref{J2-2} proves the first part of the lemma.

(ii)
From \eqref{idea}, we  obtain
\begin{align}
\mathcal R_2(\alpha,\beta,\tau)=&\int_0^\tau\Big(\fe^{is(\alpha+\beta)}-\fe^{is\alpha}
           -i\beta \fe^{is\beta} \mathcal{M}_\tau(\fe^{i\alpha\cdot})\Big)ds.
\end{align}
First, we decompose
\begin{align*}
\int_0^\tau\Big(\fe^{is(\alpha+\beta)}-\fe^{is\alpha}\Big)ds
&=\int_0^\tau\Big(\fe^{is(\alpha+\beta)}-\fe^{is\alpha}-is\beta\fe^{is\alpha} \Big)ds+i\beta \int_0^\tau s \fe^{is\alpha}\,ds\\
&=\int_0^\tau\fe^{is\alpha}\big(\fe^{is\beta}-1-is\beta\big)ds+i\beta \tau\mathcal{M}_\tau(\fe^{i\alpha\cdot})
\end{align*}
and thus get
\begin{align}\label{a2}
   \mathcal R_2(\alpha,\beta,\tau)
    =&\int_0^\tau\fe^{is\alpha}\big(\fe^{is\beta}-1-is\beta\big)\,ds
      +i\beta\int_0^\tau(1-\fe^{is\beta})\,ds\cdot \mathcal{M}_\tau(\fe^{i\alpha\cdot}).
\end{align}
Note that
$$
\Big|\fe^{is\alpha}\big(\fe^{is\beta}-1-is\beta\big)\Big|\lesssim s^2 |\beta|^2,\quad
\Big|1-\fe^{is\beta}\Big|\lesssim s|\beta|,\quad
|\mathcal{M}_\tau(\fe^{i\alpha\cdot})|\lesssim \tau.
$$
Therefore, \eqref{a2} can be controlled by $C\tau^3|\beta|^2$.
\end{proof}

\section{Construction of the method}\label{sec:construction}
Now we derive a second-order numerical method for \eqref{model}. Since the employed approximations form the basis of our error analysis, we present some details of the construction. For an alternative derivation of this method, we refer to \cite{bs}.

Let  $\tau>0$ be the time step size and $t_n=n\tau$, $n\ge 0$ the temporal grid points. First, by employing the twisted variable $v=\fe^{-it\Delta}u$ and Duhamel's formula, we get
\begin{equation}\label{v-eq-s}
 v(t_n+\sigma)=v(t_n)+i\int_0^\sigma\fe^{-i(t_n+\rho)\Delta}
  \Big(\big|\fe^{i(t_n+\rho)\Delta}v(t_n+\rho)\big|^2\fe^{i(t_n+\rho)\Delta}v(t_n+\rho)\Big)d\rho.
\end{equation}
Then, freezing the nonlinear interaction by approximating  $\fe^{i(t_n+\rho)\Delta}\approx \fe^{i(t_n+\sigma)\Delta}$ and $v(t_n+\rho)\approx v(t_n)$, we get
\begin{align}\label{r3}
  v(t_n+\sigma)= v(t_n)+i\sigma \fe^{-i(t_n+\sigma)\Delta}
 \Big(\big|\fe^{i(t_n+\sigma)\Delta}v(t_n)\big|^2\fe^{i(t_n+\sigma)\Delta}v(t_n)\Big)
 +\mathcal R_3^n(v,\sigma).
 \end{align}
The remainder term $\mathcal R_3^n(v,\sigma)$ satisfies the following estimate.

\begin{lemma}\label{lem:r3}
Let $\gamma>\frac{d}{2}$, $\sigma\in [0,\tau]$ and $v\in L^\infty((0,T);H^{\gamma+2})$. Then,
$$
\big\|\mathcal R_3^n(v,\sigma)\big\|_{H^{\gamma}}
\lesssim \tau^2\big(\|v\|_{L^\infty((0,T);H^{\gamma+2})}+ \|v\|_{L^\infty((0,T);H^{\gamma+2})}^3\big).
$$
\end{lemma}
We postpone the proof of the lemma to Section \ref{sec: Remainder-term}.

Next, we derive a second-order expansion of Duhamel's formula
\begin{equation}\label{v eq}
  v(t_n+\tau)=v(t_n)+i\int_0^\tau\fe^{-i(t_n+\sigma)\Delta}
  \Big(\big|\fe^{i(t_n+\sigma)\Delta}v(t_n+\sigma)\big|^2\fe^{i(t_n+\sigma)\Delta}v(t_n+\sigma)\Big)d\sigma.
\end{equation}
Replacing $v(t_n+\sigma)$ by \eqref{r3}, we infer that
 \begin{align}\label{r4}
 v(t_{n+1}) = v(t_n)+I_1(t_n)+I_2(t_n)+\mathcal R_4^n(v),
\end{align}
where
\begin{align}
I_1(t_n)&=i\int_0^\tau\fe^{-i(t_n+s)\Delta}
            \Big(\big|\fe^{i(t_n+s)\Delta}v(t_n)\big|^2
                \fe^{i(t_n+s)\Delta}v(t_n)\Big)ds,\notag\\
I_2(t_n)&=-\int_0^\tau s\fe^{-i(t_n+s)\Delta}
             \Big(\big|\fe^{i(t_n+s)\Delta}v(t_n)\big|^4
                \fe^{i(t_n+s)\Delta}v(t_n)\Big)ds.\label{def-I2}
\end{align}
The remainder term $\mathcal R_4^n(v)$ can be bounded as stated in the next lemma. Again, the proof of this lemma is postponed to Section~\ref{sec: Remainder-term}.

\begin{lemma}\label{lem:r4}
Let $\gamma >\frac{d}{2}$ and $0<\tau\leq1$. Then, for $v\in L^\infty((0,T);H^{\gamma+2})$,
$$
\big\|\mathcal{R}_4^n(v)\big\|_{H^\gamma}\le C\tau^3,
$$
where the constant $C$ only depends on $\|v\|_{L^\infty((0,T);H^{\gamma+2})}$.
\end{lemma}

Due to the complexity of the phase functions
$$
\phi_3=|\bxi|^2+|\bxi_1|^2-|\bxi_2|^2-|\bxi_3|^2,\quad
\phi_5=|\bxi|^2+|\bxi_1|^2+|\bxi_2|^2-|\bxi_3|^2-|\bxi_4|^2-|\bxi_5|^2,
$$
we note that the terms  in $I_1$ and $I_2$  can not be  easily expressed in physical space.

Therefore, we consider $I_1$ first in Fourier space and write
\begin{align*}
\widehat{I_1}(t_n,\bxi)=i\int_0^\tau\!\! \int_{\bxi=\bxi_1+\bxi_2+\bxi_3}\fe^{i(t_n+s)\phi_3}
           \hat{\bar v}(t_n,\bxi_1)\hat v(t_n,\bxi_2)\hat v(t_n,\bxi_3)\, (d\bxi_1)(d\bxi_2)\,ds.
\end{align*}
The main problem concerns the handling of the phase $\fe^{is\phi_3}$. Defining
\begin{equation*}
\alpha=2|\bxi_1|^2,\quad \beta=2\bxi_1\cdot\bxi_2+2\bxi_1\cdot\bxi_3+2\bxi_2\cdot\bxi_3
\end{equation*}
allows us to write
\begin{align*}
\fe^{is\phi_3}=\fe^{is\alpha+is\beta}.
\end{align*}
Applying the formulas presented in Lemma \ref{lem:idea-1}, we get
\begin{align}{\label{i1}}
\widehat{I_1}(t_n,\bxi)&=i\tau\int_{\bxi=\bxi_1+\bxi_2+\bxi_3} \varphi(i\tau \alpha)
          \fe^{it_n\phi_3} \hat{\bar v}(t_n,\bxi_1)\hat v(t_n,\bxi_2)\hat v(t_n,\bxi_3)\, (d\bxi_1)(d\bxi_2)\notag\\
           &\quad-i\tau \int_{\bxi=\bxi_1+\bxi_2+\bxi_3}\big(\fe^{i\tau\beta}-1\big)\psi(i\tau \alpha) \fe^{it_n\phi_3}\hat{\bar v}(t_n,\bxi_1)\hat v(t_n,\bxi_2)\hat v(t_n,\bxi_3)\, (d\bxi_1)(d\bxi_2)\notag\\
          & \quad+\widehat{\mathcal R_5^n}(v)(\bxi),
\end{align}
where the remainder term $\mathcal R_5^n(v)$ obeys the bound given in the following lemma. Its proof will be postponed to Section~\ref{sec: Remainder-term}.

\begin{lemma}\label{lem:r5}
Let $\gamma>\frac{d}{2}$ and $v\in L^\infty((0,T);H^{\gamma+2})$.   Then,
$$
\big\|\mathcal R_5^n(v)\big\|_{H^{\gamma}}
\lesssim \tau^3 \big\|v\big\|_{L^\infty((0,T);H^{\gamma+2})}^3.
$$
\end{lemma}

Using $\beta=\phi_3-\alpha$, we transform \eqref{i1} back to physical space to get
\begin{align}
I_1&=i\tau \fe^{-it_n\Delta}\Big\{ \big[ \varphi(-2i\tau \Delta)\fe^{-it_n\Delta}\bar v(t_n)\big]\cdot\Big(\fe^{it_n\Delta}v(t_n)\Big)^2\Big\}\notag\\
&\quad-i\tau  \fe^{-it_{n+1}\Delta}\Big\{\big[ \psi(-2i\tau \Delta)\fe^{-it_{n-1}\Delta}\bar v(t_n)\big]\cdot\Big(\fe^{it_{n+1}\Delta}v(t_n)\Big)^2\Big\}\label{def-I1-pointwise}\\
&\quad+i\tau \fe^{-it_n\Delta}\Big\{ \big[ \psi(-2i\tau \Delta)\fe^{-it_n\Delta}\bar v(t_n)\big]\cdot\Big(\fe^{it_n\Delta}v(t_n)\Big)^2\Big\}+\mathcal R_5^n(v)\notag.
\end{align}
The term $I_2$ is of higher order in $\tau$. Therefore, it is sufficient to freeze the linear flow and approximate the term as
 \begin{align}
I_2(t_n)&=-\int_0^\tau s\fe^{-it_n\Delta}
             \Big(\big|\fe^{it_n\Delta}v(t_n)\big|^4
                \fe^{it_n\Delta}v(t_n)\Big)ds+\mathcal R_6^n(v)\label{def-R6}\\
            &=
               -\frac12 \tau^2\fe^{-it_n\Delta}
             \Big(\big|\fe^{it_n\Delta}v(t_n)\big|^4
                \fe^{it_n\Delta}v(t_n)\Big)+\mathcal R_6^n(v),\label{def-I2-pointwise}
\end{align}
where the remainder term $\mathcal R_6^n(v)$ obeys the bound given in the following lemma. Again, its proof will be postponed to Section \ref{sec: Remainder-term}.

\begin{lemma}\label{lem:r6}
Let $\gamma>\frac{d}{2}$ and $v\in L^\infty((0,T);H^{\gamma+2})$. Then
$$
\big\|\mathcal R_6^n(v)\big\|_{H^{\gamma}}
\lesssim \tau^3 \|v\|_{L^\infty((0,T);H^{\gamma+2})}^5.
$$
\end{lemma}

Now combining  \eqref{r4}, \eqref{def-I1-pointwise} and \eqref{def-I2-pointwise}, we have that
\begin{align}\label{v-n+1}
v(t_{n+1})= \Phi^n\big(v(t_n)\big)+\mathcal R_4^n(v)+\mathcal R_5^n(v)+\mathcal R_6^n(v),
\end{align}
where the operator $\Phi^n$ is defined by
\begin{equation}\label{def-Phi-f}
\begin{aligned}
\Phi^n\big(f\big)
=f&+i\tau \fe^{-it_n\Delta}\Big\{ \Big( \varphi(-2i\tau \Delta)\fe^{-it_n\Delta}\bar f\Big)\cdot\Big(\fe^{it_n\Delta}f\Big)^2\Big\}\\
&-i\tau  \fe^{-it_{n+1}\Delta}\Big\{\Big( \psi(-2i\tau \Delta)\fe^{-it_{n-1}\Delta}\bar f\Big)\cdot\Big(\fe^{it_{n+1}\Delta}f\Big)^2\Big\}\\
&+i\tau \fe^{-it_n\Delta}\Big\{ \Big( \psi(-2i\tau \Delta)\fe^{-it_n\Delta}\bar f\Big)\cdot\Big(\fe^{it_n\Delta}f\Big)^2\Big\}\\
&-\tfrac12 \tau^2\fe^{-it_n\Delta} \Big(\big|\fe^{it_n\Delta}f\big|^4 \fe^{it_n\Delta}f\Big).
\end{aligned}
\end{equation}

Our second order low-regularity integrator is obtained by dropping the remainder terms $\mathcal R_4^n,\mathcal R_5^n,\mathcal R_6^n$ in \eqref{v-n+1}. The method for the twisted variable is summarized as follows: let $v^0=u_0$ and
\begin{align} \label{NLRI2-v}
 v^{n+1}=&\,\Phi^n\big(v^n\big) \quad \mbox{ for } n\geq 0.
\end{align}
Finally, setting $u^n=\fe^{it_n\Delta}v^n$, we obtain the announced numerical scheme \eqref{LRI-discrete-1} for the NLS equation \eqref{model}.

\subsection{Estimates of the remainder terms}\label{sec: Remainder-term}
Now we prove Lemmas \ref{lem:r3} to \ref{lem:r6}.

\begin{proof}[Proof of Lemma \ref{lem:r3}]
By \eqref{r3}, we have that
\begin{align*}
\mathcal R_3^n&(v,s)
=i\int_0^s \Big(\fe^{-i(t_n+\sigma)\Delta}-\fe^{-i(t_n+s)\Delta}\Big)
  \Big(\big|\fe^{i(t_n+\sigma)\Delta}v(t_n+\sigma)\big|^2\fe^{i(t_n+\sigma)\Delta}v(t_n+\sigma)\Big)\,d\sigma\\
&+ i\!\int_0^s \fe^{-i(t_n+s)\Delta}
  \Big(\big|\fe^{i(t_n+\sigma)\Delta}v(t_n+\sigma)\big|^2-\big|\fe^{i(t_n+s)\Delta}v(t_n+s)\big|^2\Big)\fe^{i(t_n+\sigma)\Delta}v(t_n+\sigma)\,d\sigma\\
&+i\!\int_0^s\! \fe^{-i(t_n+s)\Delta}
  \Big(\big|\fe^{i(t_n+s)\Delta}v(t_n+s)\big|^2\big(\fe^{i(t_n+\sigma)\Delta}v(t_n+\sigma)-\fe^{i(t_n+s)\Delta}v(t_n)\big)\Big)\,d\sigma.
\end{align*}
Note that from \eqref{v-eq-s}, Lemma \ref{lem:kato-Ponce}\;(i) and Sobolev embedding, we get
$$
\sup_{0\le\sigma\le\tau}\big\| v(t_n+\sigma)-v(t_n)
\big\|_{H^{\gamma}}
\lesssim \tau
\|v\|_{L^\infty((0,T);H^{\gamma})}^3.
$$
Moreover, for any $f\in H^\gamma$,
\begin{align}\label{est:t-s}
\left\|\Big(\fe^{-i(t_n+\sigma)\Delta}-\fe^{-i(t_n+s)\Delta}\Big)f\right\|_{H^\gamma}
\lesssim
|\sigma-s| \|f\|_{H^{\gamma+2}}.
\end{align}
Applying these two estimates, we obtain
$$
\big\|\fe^{i(t_n+\sigma)\Delta}v(t_n+\sigma)-\fe^{i(t_n+s)\Delta}v(t_n)\big\|_{H^\gamma}
\lesssim \tau \big(\|v\|_{L^\infty((0,T);H^{\gamma+2})}+\|v\|_{L^\infty((0,T);H^{\gamma+2})}^3\big)
$$
and thus
\begin{align*}
\big\|\mathcal R_3^n(v,s)\big\|_{H^\gamma} \lesssim
\tau^2 \big(\|v\|_{L^\infty((0,T);H^{\gamma+2})}+\|v\|_{L^\infty((0,T);H^{\gamma+2})}^3\big).
\end{align*}
This is the desired result.
\end{proof}

\begin{proof}[Proof of Lemma \ref{lem:r4}]
Inserting \eqref{r3} with $\sigma=\rho$ in \eqref{v-eq-s} and using \eqref{r4} shows that the remainder $\mathcal{R}_4^n(v)$ consists of terms of the form
$$
i\int_0^\tau\fe^{-i(t_n+s)\Delta}
            \Big(\fe^{i(t_n+s)\Delta}\mathcal W_j\cdot\fe^{-i(t_n+s)\Delta}\overline{\mathcal W_k}\cdot\fe^{i(t_n+s)\Delta}\mathcal W_\ell\Big)ds,\qquad j+k+\ell\ge 5,
$$
where
\begin{align*}
\mathcal W_1 = v(t_n),\quad
\mathcal W_2 = is\fe^{-i(t_n+s)\Delta}\Big(\big|\fe^{i(t_n+s)\Delta}v(t_n)\big|^2\fe^{i(t_n+s)\Delta}v(t_n)\Big), \quad
\mathcal W_3 = \mathcal R_3^n(v,s).
\end{align*}
By Lemma \ref{lem:r3} and Lemma \ref{lem:kato-Ponce}\;(i), we thus get
 \begin{align*}
\big\|\mathcal{R}_4^n(v)\big\|_{H^\gamma}
\lesssim  C\Big(\big\|v(t_n) \big\|_{L^\infty((0,t);H^{\gamma+2})}\Big)\tau^3.
\end{align*}
This finishes the proof of the lemma.
\end{proof}

\begin{proof}[Proof of Lemma \ref{lem:r5}]
Without loss of generality, we may assume that $\hat v(t_n)$ and $\hat{\bar v}(t_n)$ are positive (otherwise one may replace them by their absolute values).

From Lemma \ref{lem:idea-1}, we have
\begin{align*}
\widehat{\mathcal R_5^n}(v)(\bxi) =
\int_{\bxi=\bxi_1+\bxi_2+\bxi_3}\mathcal R_2(\alpha,\beta,\tau)\,
 \fe^{it_n\phi_3} \hat {\bar v}(t_n,\bxi_1)\hat v(t_n,\bxi_2)\hat v(t_n,\bxi_3)\,(d\bxi_1)(d\bxi_2)
\end{align*}
and further
\begin{align*}
\big|\widehat{\mathcal R_5^n}(v)(\bxi)\big|
\lesssim
\tau^3\int_{\bxi=\bxi_1+\bxi_2+\bxi_3} \beta^2\,\hat {\bar v}(t_n,\bxi_1)\hat v(t_n,\bxi_2)\hat v(t_n,\bxi_3)\,(d\bxi_1)(d\bxi_2).
\end{align*}
By symmetry, we may assume that $|\bxi_1|\ge |\bxi_2|\ge |\bxi_3|$. This yields
\begin{align*}
\langle\bxi\rangle^{\gamma}\beta^2
&\lesssim
 \langle\bxi\rangle^{\gamma} \big(|\bxi_1|^2|\bxi_2|^2+|\bxi_1|^2|\bxi_3|^2+|\bxi_2|^2|\bxi_3|^2\big)\\
&\lesssim |\bxi_1|^{2+\gamma}|\bxi_2|^2.
\end{align*}
Using this estimate, we get
\begin{align*}
\langle\bxi\rangle^{\gamma}\big|
\widehat{\mathcal R_5^n}&(v)(\bxi)\big|\\
 &\lesssim\tau^3\int_{\bxi=\bxi_1+\bxi_2+\bxi_3,|\bxi_1|\ge |\bxi_2|\ge |\bxi_3|}|\bxi_1|^{2+\gamma}|\bxi_2|^2\hat {\bar v}(t_n,\bxi_1)\hat v(t_n,\bxi_2)\hat v(t_n,\bxi_3)\,(d\bxi_1)(d\bxi_2)\\
&\lesssim  \tau^3\mathcal F\Big(\big(\minus\Delta\big)^{1+\gamma/2}\bar v\cdot \big(\minus\Delta\big)v\cdot v\Big)(t_n,\bxi).
\end{align*}
Therefore, by Plancherel's identity and Lemma \ref{lem:kato-Ponce}\;(ii) with $\delta=0$, we obtain that for any $\gamma_1>\frac{d}{2}$,
\begin{align*}
\left\|\mathcal R_5^n(v)\right\|_{H^{\gamma}}
&\lesssim
\tau^3   \left\|\big(\minus\Delta\big)^{1+\gamma/2}\bar v\cdot \big(\minus\Delta\big)v\cdot v\right\|_{L^\infty((0,T);L^2)}\\
&\lesssim  \tau^3   \big\|v\big\|_{L^\infty((0,T);H^{\gamma+2})}\big\|v\big\|_{L^\infty((0,T);H^{\gamma_1+2})}
       \big\|v\big\|_{L^\infty((0,T);H^{\gamma_1})}.
\end{align*}
Since $\gamma>\frac{d}{2}$, choosing $\gamma_1=\gamma$, we get the desired result.
\end{proof}

\begin{proof}[Proof of Lemma \ref{lem:r6}]
By \eqref{def-I2} and \eqref{def-R6}, we have that
\begin{align*}
\mathcal R_6^n=
&-\int_0^\tau s\Big(\fe^{-i(t_n+s)\Delta}-\fe^{-it_n\Delta}\Big)
             \Big(\big|\fe^{i(t_n+s)\Delta}v(t_n)\big|^4
                \fe^{i(t_n+s)\Delta}v(t_n)\Big)ds\\
&-\int_0^\tau s\,\fe^{-it_n\Delta}
             \Bigl(\Big[\big|\fe^{i(t_n+s)\Delta}v(t_n)\big|^4 - \big|\fe^{it_n\Delta}v(t_n)\big|^4\Big]
                \fe^{i(t_n+s)\Delta}v(t_n)\Big)ds\\
&-\int_0^\tau s\,\fe^{-it_n\Delta}
             \Big(\big|\fe^{it_n\Delta}v(t_n)\big|^4\cdot
                \big(\fe^{-i(t_n+s)\Delta}-\fe^{-it_n\Delta}\big)v(t_n)\Big)ds.
\end{align*}
Then, the claimed result follows directly from \eqref{est:t-s} and Lemma \ref{lem:kato-Ponce}\;(i).
\end{proof}

\section{Proof of Theorem \ref{thm:convergence}}\label{sec:ord-proof}
Taking the difference between the numerical scheme \eqref{NLRI2-v} and the exact solution gives
\begin{align*}
v^{n+1}-v(t_{n+1})&=\Phi^n\big(v(t_n)\big)-v(t_{n+1})+ \Phi^n\big(v^n\big)-\Phi^n\big(v(t_n)\big)\\
         &= \mathcal{L}^n+\Phi^n\big(v^n\big)-\Phi^n\big(v(t_n)\big),
\end{align*}
where $\mathcal{L}^n=\Phi^n\big(v(t_n)\big)-v(t_{n+1})$ is the local error.

\subsection{Local error}
\allowdisplaybreaks
The following bound on the local error holds.
\begin{lemma}\label{lem:local-error} Let $\gamma >\frac{d}{2}$ and $0<\tau\leq1$. Then,
$$
\big\|\mathcal{L}^n\big\|_{H^\gamma}\le C\tau^3,
$$
 where the constant  $C$ only depends on $\|v\|_{L^\infty((0,T);H^{\gamma+2})}$.
\end{lemma}
\begin{proof}
By \eqref{v-n+1}, we get that
$$
\mathcal{L}^n=-\mathcal{R}_4^n(v)-\mathcal R_5^n(v)-\mathcal R_6^n(v).
$$
Thus, the desired estimate follows from Lemmas \ref{lem:r4}, \ref{lem:r5}, and \ref{lem:r6}.
\end{proof}

\subsection{Stability}
The main result in this subsection is the following stability estimate.
\begin{lemma}\label{lem:stability} Let $\gamma >\frac d2$. Then,
$$
\big\|\Phi^n(v^n)-\Phi^n(v(t_n))\big\|_{H^\gamma}\le (1+ C\tau)\big\|v^n-v(t_n)\big\|_{H^\gamma}+C\tau\big\|v^n-v(t_n)\big\|_{H^\gamma}^5,
$$
where the constant $C$ only depends on $\|v\|_{L^\infty((0,T);H^{\gamma})}$.
\end{lemma}
\begin{proof}
For short, we denote $g_n=v^n-v(t_n)$. Then, using \eqref{def-Phi-f}, we have
\begin{align*}
\Phi^n(v^n)&-\Phi^n(v(t_n))=g_n+\sum\limits_{j=1}^4\big(\Phi^n_j(v^n)-\Phi^n_j(v(t_n))\big),
\end{align*}
where
\begin{align*}
\Phi^n_1(f)&=i\tau \fe^{-it_n\Delta}\Big\{ \Big( \varphi(-2i\tau \Delta)\fe^{-it_n\Delta}\bar f\Big)\cdot\Big(\fe^{it_n\Delta}f\Big)^2\Big\}\\
\Phi^n_2(f)&=-i\tau  \fe^{-it_{n+1}\Delta}\Big\{\Big( \psi(-2i\tau \Delta)\fe^{-it_{n-1}\Delta}\bar f\Big)\cdot\Big(\fe^{it_{n+1}\Delta}f\Big)^2\Big\}\\
\Phi^n_3(f)&=i\tau \fe^{-it_n\Delta}\Big\{ \Big( \psi(-2i\tau \Delta)\fe^{-it_n\Delta}\bar f\Big)\cdot\Big(\fe^{it_n\Delta}f\Big)^2\Big\}\notag\\
\Phi^n_4(f)&=-\frac12 \tau^2\fe^{-it_n\Delta}\Big(\big|\fe^{it_n\Delta}f\big|^4\fe^{it_n\Delta}f\Big).
\end{align*}

Note that  by the definition of $\varphi$ and $\psi$ in \eqref{11}, we have that
\begin{align*}
\big\|\varphi(-2i\tau \Delta)f\big\|_{H^\gamma}
\lesssim \|f\|_{H^\gamma},\qquad
\big\|\psi(-2i\tau \Delta)f\big\|_{H^\gamma}
\lesssim \|f\|_{H^\gamma}.
\end{align*}
Hence, by Lemma \ref{lem:kato-Ponce}\,(i),
\begin{align*}
\big\|\Phi^n_1(v^n)-\Phi^n_1(v(t_n))\big\|_{H^\gamma}
\leq & C\tau
\Big(\big\|g_n\big\|_{H^{\gamma}}+\big\|g_n\big\|_{H^{\gamma}}^3\Big),
\end{align*}
where
$C$ only depends on $\big\|v\big\|_{L^\infty((0,T); H^{\gamma})}$.

Similarly, we get that
\begin{align}
\sum_{j=2}^4\big\|\Phi^n_j(v^n)-\Phi^n_j(v(t_n))\big\|_{H^\gamma}
\leq C\tau
\Big(\big\|g_n\big\|_{H^{\gamma}}+\big\|g_n\big\|_{H^{\gamma}}^5\Big).\label{est:Phi-2}
\end{align}
Combining the above estimates, we finally obtain
\begin{align*}
\big\|\Phi^n\big(v(t_n)\big)-\Phi^n\big(v^n\big)\big\|_{H^\gamma}
\le \|g_n\|_{H^\gamma}+C\tau \Big(\big\|g_n\big\|_{H^{\gamma}}+\big\|g_n\big\|_{H^{\gamma}}^5\Big),
\end{align*}
which is the desired result.
\end{proof}

\subsection{Proof of Theorem \ref{thm:convergence}}
Now, combining the local error estimate with the stability result, we give the proof of Theorem \ref{thm:convergence}. From Lemma \ref{lem:local-error} and Lemma \ref{lem:stability}, we infer that there exits a constant $C$ depending only on $\|v\|_{L^\infty((0,T);H^{\gamma+2})}$, such that for $0<\tau\leq1$, we have
\begin{align*}
\big\|v(t_{n+1})-v^{n+1}\big\|_{H^\gamma}
\le C\tau^3+(1+C\tau)\big\|v(t_{n})-v^{n}\big\|_{H^\gamma}+C\tau\big\|v(t_{n})-v^{n}\big\|_{H^\gamma}^5, \quad  n\geq 0.
\end{align*}
By recursion, we get from this the bound
$$
\big\|v(t_{n+1})-v^{n+1}\big\|_{H^\gamma}
\le C\tau \sum_{j=0}^n (1+C\tau)^j\Big[\big\|v(t_{n-j})-v^{n-j}\big\|_{H^\gamma}^5+C\tau^2\Big].
$$
From this estimate we infer that there exist positive constants $\tau_0$ and $C$, such that for any $\tau\in[0,\tau_0]$,
\begin{align*}
\big\|v(t_{n+1})-v^{n+1}\big\|_{H^\gamma}
\le C\tau^3\sum\limits_{j=0}^n(1+C\tau)^j\le C \tau^2,\quad n\geq 0.
\end{align*}
Note that the constants $\tau_0$ and $C$ only depend on $T$ and $\|u\|_{L^\infty((0,T);H^{\gamma+2})}$.
This proves Theorem \ref{thm:convergence}.\qed

\section{Numerical experiments} \label{sec:numerical}
In this section we carry out some numerical experiments to illustrate our convergence result in two space dimensions. For this purpose,
we consider the nonlinear Schr\"{o}dinger equation \eqref{model} with initial data
\begin{align}\label{initial}
	u_0(x_1,x_2)=\sum_{(k,\ell)\in\Z^2}\big(1+\sqrt{k^2+\ell^2}\big)^{-\frac12-\gamma-\varepsilon}(1+i) \,\fe^{i(kx_1+\ell x_2)},\qquad\varepsilon>0,
\end{align}
where $\gamma$ is used to set the regularity of the data. This choice guarantees that $u^0\in H^\gamma(\mathbb T^2)$. In the experiment, we set $\varepsilon=0$.

We choose $N=2^7$, i.e.~$2^{14}$ grid points, and measure the temporal discretisation error $w = u(t_n,\cdot)-u^n_{\tau,N}$ in the discrete $H^\gamma$-norm
$$
\|w\|_{H^\gamma_N} = \left\|\big(1+ \sqrt{-\Delta}\big)^\gamma w\right\|_{L^2_N},
$$
where
$$
\|w\|_{L^2_N}^2=\frac{4\pi^2}{N^2}\sum_{j,m=0}^{N-1}|w(x_1^j, x_2^m)|^2, \quad x_1^j=\frac{2\pi j}{N},\ x_2^m=\frac{2\pi m}{N}.
$$
Our results for initial data $u_0\in H^{\gamma+2}(\mathbb T^2)$ are presented in Fig.~1. We choose the three different values $\gamma=1$, $1.5$, $2$ to illustrate the convergence rate of our scheme \eqref{LRI-discrete-1}. As expected, the slopes of the error curves are 2 whenever $\gamma$ is bigger than 1. The slope of the curve for $\gamma=1$ is slightly less regular. This is also expected because the value $\gamma=1$ is the limit case in two space dimensions. Thus, the results  agree well with the corresponding results of the theoretical analysis, given in Theorem \ref{thm:convergence}.

\begin{figure}[!htb]
\includegraphics[height=7.5cm, keepaspectratio]{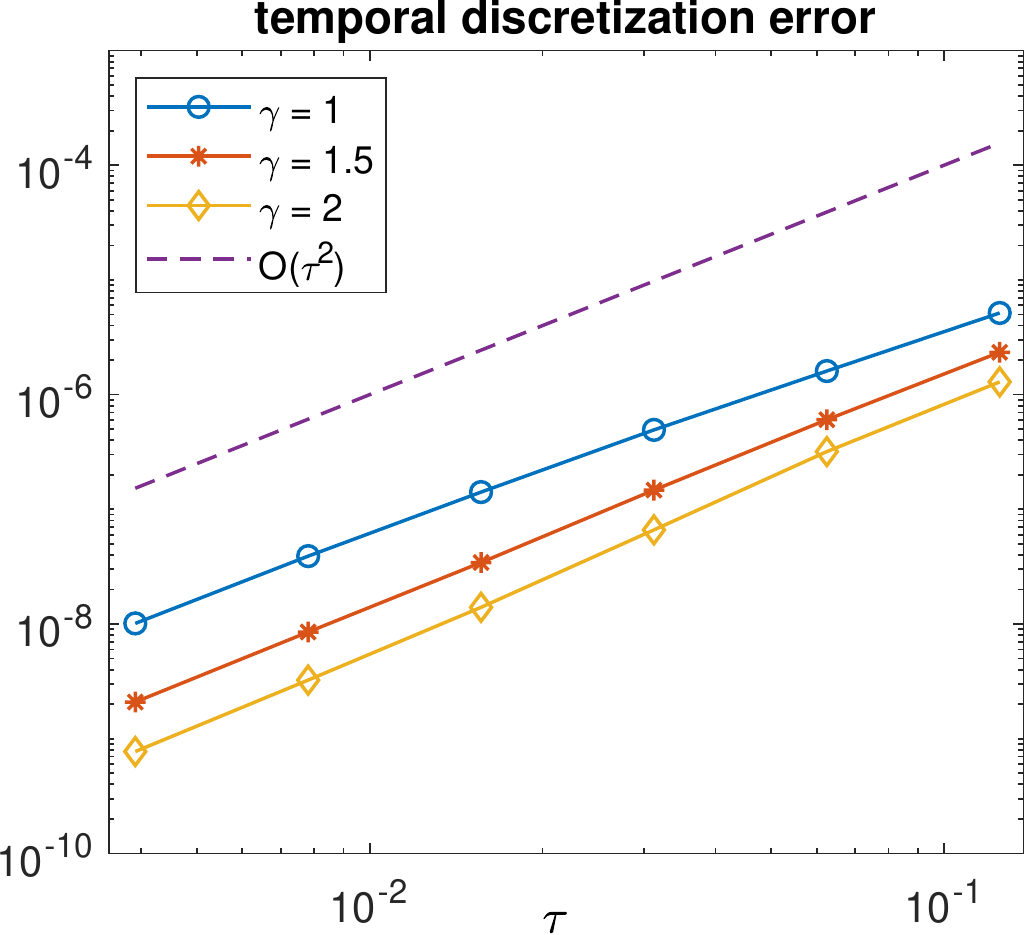}
\caption{ Temporal discretisation error in $H^\gamma$ for initial data in $H^{\gamma+2}$ for various values of $\gamma$. The errors are measured at $T=1$ for various step sizes $\tau$. The dashed line has slope 2.}\label{fig:NLRI}
\end{figure}

\section*{Declarations}

\noindent
{\bf Availability of data and materials}\\
Not applicable.

\medskip
\noindent
{\bf Competing interests}\\
The authors declare that they have no competing interests.

\newpage
\medskip
\noindent
{\bf Funding}\\
Y.W.~was partially supported by the NSFC grants 12171356 and 11771325.
F.Y.~acknowledges financial support by the China Scholarship Council.

\medskip
\noindent
{\bf Authors’ contributions}\\
The authors declare that the study was realized in collaboration with equal responsibility. All authors read and approved the final manuscript.

\medskip
\noindent
{\bf Acknowledgements}\\
Parts of the research were carried out during a research visit of F.Y.~at the University of Innsbruck.

\bibliographystyle{model1-num-names}

\end{document}